\documentclass[12pt,fleqn]{article}
\usepackage[cp866]{inputenc}
\usepackage[english]{babel}
\usepackage{makeidx}
\usepackage{latexsym,amsfonts,amssymb,amsmath,longtable,amsthm}
\input amssym.def
\usepackage{latexsym}
\usepackage{graphicx}
\usepackage[all]{xy}
\usepackage{amsfonts}
\usepackage{amsmath}
\usepackage{mathrsfs}
\usepackage{color}
\usepackage{ulem}
\usepackage[table]{xcolor}

\newcommand{\inter}{\mathop{\rm Int}}

\newcommand{\cl}{\mathop{\rm Cl}}

\newcommand{\Proj}{\mathop{\rm Proj_{rot}}}

\begin{document}
\title{On Properties of the Intrinsic
Geometry of Submanifolds in a Riemannian Manifold\footnote{{\it Mathematical Subject classification\/} (2010).
53C45; {\it Key words}: intrinsic metric, induced boundary metric,
triangle inequality, geodesics.}}

\author{ Anatoly P.
Kopylov\footnote{Sobolev Institute of Mathematics, pr. Akad. Koptyuga 4, 630090, Novosibirsk, Russia
and Novosibirsk State University, ul. Pirogova 2, 630090,
Novosibirsk, Russia; E-mail: apkopylov@yahoo.com}, Mikhail V.
Korobkov\footnote{Sobolev Institute of Mathematics, pr. Akad. Koptyuga 4, 630090, Novosibirsk, Russia
and Novosibirsk State University, ul. Pirogova 2, 630090,
Novosibirsk, Russia; E-mail: korob@math.nsc.ru}}

\maketitle

\begin{abstract} Assume that
$(X,g)$ is an $n$-dimensional smooth connected Riemannian manifold
without boundary and $Y$ is an $n$-dimensional compact connected
$C^0$-submanifold in~$X$ with nonempty boundary $\partial Y$ ($n \ge 2$).
We consider the metric function $\rho_Y(x,y)$
generated by the intrinsic metric of the interior $\inter Y$
of~$Y$ in the following natural way: $\rho_Y(x,y) =
\liminf\limits_{x' \to x,\,y' \to x;\,\,x',y' \in \inter Y}
\{\inf[l(\gamma_{x',y',\inter Y})]\}$, where
$\inf[l(\gamma_{x',y',\inter Y})]$ is the infimum of the lengths
of smooth paths joining $x'$ and $y'$ in the interior $\inter Y$
of~$Y$. We study conditions under which $\rho_Y$ is
a metric and also the question about the existence of geodesics in
the metric $\rho_Y$ and its relationship with the classical
intrinsic metric of the hypersurface~$\partial Y$.
\end{abstract}

Let $(X,g)$ be an $n$-dimensional smooth connected Riemannian
manifold without boundary and let $Y$ be an~$n$-dimensional
compact connected $C^0$-submanifold in~$X$ with nonempty boundary
$\partial Y$ ($n \ge 2$). A classical object of investigations
(see, for example,~\cite{al1}) 
is given by the intrinsic metric
$\rho_{\partial Y}$ on the hypersurface $\partial Y$ defined for
$x,y \in \partial Y$ as the infimum of the lengths of curves $\nu
\subset \partial Y$ joining $x$ and $y$. In the recent decades, an
alternative approach arose in the rigidity theory for submanifolds
of Riemannian manifolds (see, for instance, the recent
articles~\cite{kop1}, 
\cite{kop2,kor1}, 
which also contain a historical survey of
works on the topic). In accordance with this approach, the metric
on $\partial Y$ is induced by the intrinsic metric of the interior
$\inter Y$ of the submanifold $Y$. Namely, suppose that $Y$
satisfies the following condition\footnote{Easy examples show that
if $X$ is an $n$-dimensional connected smooth Riemannian manifold
without boundary then an~$n$-dimensional compact connected
$C^0$-submanifold in $X$ with nonempty boundary may fail to
satisfy condition~$(i)$. For $n = 2$, we have the following
counterexample: Let $(X,g)$ be the space $\mathbb R^2$ equipped
with the~Euclidean metric and let $Y$ be a~closed Jordan domain in
$\mathbb R^2$ whose boundary is the union of the singleton $\{0\}$
consisting of the origin~$0$, the segment $\{(1-t)(e_1 + 2e_2) +
t(e_1 + e_2) : 0 \le t \le 1\}$, and the segments of the following
four types:\pagebreak
$$
\biggl\{\frac{(1-t)(e_1 + e_2)}{n} + \frac{te_1}{n + 1}\, :\,
0 \le t \le 1\biggr\}  \quad (n = 1,2,\dots);
$$
$$
\biggl\{\frac{e_1 + (1-t)e_2}{n}\, :\, 0 \le t \le 1\biggr\}
\quad (n = 2,3,\dots);
$$
$$
\biggl\{\frac{(1-t)(e_1 + 2e_2)}{n} + \frac{2t(2e_1 +e_2)}{4n
+3}\, :\, 0 \le t \le 1\biggr\} \quad (n = 1,2,\dots);
$$
$$
\biggl\{\frac{(1-t)(e_1 + 2e_2)}{n + 1} + \frac{2t(2e_1 + e_2)}{4n
+3}\, :\, 0 \le t \le 1\biggr\} \quad (n = 1,2,\dots).
$$
Here $e_1, e_2$ is the canonical basis in $\mathbb R^2$. By the
construction of $Y$, we have $\rho_Y(0,\; E) = \infty$ for every
$E \in Y\setminus\{0\}$.}:

(i) if
$x,y\in Y$,
then
\begin{equation}
\rho_Y(x,y)=\liminf_{x'\to x,\;y'\to y;\;x',\;y'\in\inter Y}
\{\inf\,[l(\gamma_{x',\;y',\;\inter Y})]\}<\infty, \label{1}
\end{equation}
where
$\inf\,[l(\gamma_{x',\;y',\;\inter Y})]$
is the infimum of the lengths
$l(\gamma_{x',\;y',\;\inter Y})$
of smooth paths
$\gamma_{x',\;y',\;\inter Y}:[0,1]\to\inter Y$
joining
$x'$
and
$y'$
in the interior
$\inter Y$
of
$Y$.

Note that the intrinsic metric of convex hypersurfaces in $\mathbb
R^n$ (i.e., a classical object) is an important particular case of
a~function $\rho_Y$. (To verify that, take as~$Y$ the~complement
of the convex hull of the hypersurface.) However, here there
appear some new phenomena. The following question is of primary
interest in our paper: Is the function $\rho_Y$ defined
by~(\ref{1}) a metric on~$Y$? If $n = 2$ then the
answer is `yes' (see Theorem~1 below) and if $n > 2$ then it is
`no' (see Theorem~2). Moreover, we prove that if $\rho_Y$ is a
metric (for an~arbitrary dimension $n \ge 2$) then any two points
$x,y \in Y$ may be joined by a shortest curve (geodesic) whose
length in the metric $\rho_Y$ coincides with $\rho_Y(x,y)$
(Theorem~3).

We will begin with the following result.

{\bf Theorem~1.}
{\it Let
$n = 2$.
Then, under condition {\rm(i),}
$\rho_Y$
is a metric on
$Y$.}

{\sc Proof.} It suffices to prove that $\rho_Y$ satisfies the
triangle inequality. Let $A$, $O$, and $D$ be three points on the
boundary of $Y$ (note that this case is basic because the other
cases are simpler). Consider $\varepsilon > 0$ and assume that
$\gamma_{A_{\varepsilon} O^1_{\varepsilon}} : [0,1] \to \inter Y$
and $\gamma_{O^2_{\varepsilon} D_{\varepsilon}} : [2,3] \to \inter
Y$ are smooth paths with the endpoints $A_{\varepsilon} =
\gamma_{A_{\varepsilon} O^1_{\varepsilon}}(0)$, $O^1_{\varepsilon}
= \gamma_{A_{\varepsilon} O^1_{\varepsilon}}(1)$ and
$D_{\varepsilon} = \gamma_{O^2_{\varepsilon} D_{\varepsilon}}(3)$,
$O^2_{\varepsilon} = \gamma_{O^2_{\varepsilon}
D_{\varepsilon}}(2)$ satisfying the conditions
$\rho_X(A_{\varepsilon}, A) \le \varepsilon$,
$\rho_X(D_{\varepsilon}, D) \le \varepsilon$,
$\rho_X(O^j_{\varepsilon}, O) \le \varepsilon$ ($j = 1;\; 2$),
$|l(\gamma_{A_{\varepsilon} O^1_{\varepsilon}}) - \rho_Y(A,O)| \le
\varepsilon$, and $|l(\gamma_{O^2_{\varepsilon} D_{\varepsilon}})
- \rho_Y(O,D)| \le \varepsilon$. Let $(U, h)$ be a chart of the
manifold $X$ such that $U$ is an open neighborhood of the point
$O$ in $X$, $h(U)$ is the unit disk $B(0, 1) = \{(x_1, x_2) \in
\mathbb R^2 : x_1^2 + x_2^2 < 1\}$ in $\mathbb R^2$, and $h(O) =
0$ ($0 = (0, 0)$ is the origin in $\mathbb R^2$); moreover, $h : U
\to h(U)$ is a diffeomorphism having the following property: there
exists a chart $(Z,\psi)$ of $Y$ with $\psi(O) = 0$, $A,D \in U
\setminus \cl_X Z$ ($\cl_X Z$ is the closure of $Z$ in the space
$(X,g)$) and $Z = \widetilde{U} \cap Y$ is the intersection of an
open neighborhood $\widetilde{U}$ ($\subset U$) of $O$ in $X$ and
$Y$ whose image $\psi(Z)$ under $\psi$ is the half-disk $B_+(0, 1)
= \{(x_1, x_2) \in B(0, 1) : x_1 \ge 0\}$. Suppose that $\sigma_r$
is an arc of the circle $\partial B(0, r)$ which is a connected
component of the set $V \cap \partial B(0, r)$, where $V = h(Z)$
and $0 < r < r^* = \min \{|h(\psi^{-1}(x_1, x_2))| : x_1^2 + x_2^2
= 1/4, \, x_1 \ge 0\}$. Among these components, there is at least
one (preserve the notation $\sigma_r$ for it) whose ends belong to
the sets $h(\psi^{-1}(\{-te_2 : 0 < t < 1\}))$ and
$h(\psi^{-1}(\{te_2 : 0 < t < 1\}))$ respectively. Otherwise, the
closure of the connected component of the set $V \cap B(0, r)$
whose boundary contains the origin would contain a point belonging
to the arc $\{e^{i\theta}/2 : |\theta| \le \pi/2\}$ (here we use
the complex notation $z = re^{i\theta}$ for points $z \in \mathbb
R^2$ ($= \mathbb C$)). But this is impossible. Therefore, the
above-mentioned arc $\sigma_r$ exists.

It is easy to check that if $\varepsilon$ is sufficiently small
then the images of the paths $h \circ
\gamma_{A_{\varepsilon}O^1_{\varepsilon}}$ and $h \circ
\gamma_{O^2_{\varepsilon}D_{\varepsilon}}$, also
intersect~$\sigma_r$, i.e., there are $t_1 \in ]0,1[$, $t_2 \in
]2,3[$ such that $\gamma_{A_\varepsilon O^1_\varepsilon}(t_1) =
x^1 \in Z$, $\gamma_{O^2_\varepsilon D_\varepsilon}(t_2) = x^2 \in
Z$ and $h(x^j) \in \sigma_r$, $j = 1,2$. Let $\widetilde\gamma_r :
[t_1,t_2] \to \sigma_r$ be a smooth parametrization of the
corresponding subarc of $\sigma_r$, i.e., $\widetilde\gamma_r(t_j)
= h(x^j)$, $j = 1,2$. Now we can define a mapping
$\gamma_\varepsilon : [0,3] \to \inter Y$ by setting
$$\gamma_\varepsilon(t)=\left\{ \aligned
\gamma_{A_{\varepsilon} O^1_{\varepsilon}}(t) , & \ \ t\in [0,t_1]; \\
h^{-1}(\widetilde\gamma_r(t)), & \ \ t\in]t_1,t_2[; \\
\gamma_{O^2_{\varepsilon}D_{\varepsilon} }(t) , & \ \ t\in
[t_2,3].
\endaligned \right.
$$
By construction, $\gamma_\varepsilon$ is a piecewise smooth path
joining the points $A_\varepsilon = \gamma_\varepsilon(0)$,
$D_\varepsilon = \gamma_\varepsilon(3)$ in $\inter Y$; moreover,
$$
l(\gamma_\varepsilon) \le l(\gamma_{A_\varepsilon O^1_\varepsilon}) +
l(\gamma_{O^2_\varepsilon D_\varepsilon}) + l(h^{-1}(\sigma_r)).
$$
By an appropriate choice of $\varepsilon > 0$, we can make $r> 0$
arbitrarily small, and since a piecewise smooth path can be
approximated by smooth paths, we have $\rho_Y(A,D) \le \rho_Y(A,O)
+ \rho_Y(O,D)$, q.e.d.

In connection with Theorem~1, there appears a natural question: Are there
analogs of this theorem for $n \ge 3$? According to the following Theorem~2,
the answer to this question is negative.

{\bf Theorem~2.} {\it If $n \ge 3$ then there exists an
$n$-dimensional compact connected $C^0$-manifold $Y \subset
\mathbb R^n$ with nonempty boundary $\partial Y$ such that
condition}~(i) ({\it where now $X = \mathbb R^n$}) {\it is
fulfilled for~$Y$ but the function $\rho_Y$ in this condition is
not a metric on $Y$.}

{\sc Proof.} It suffices to consider the case of $n = 3$. Suppose
that $A$, $O$, $D$ are points in $\mathbb R^3$, $O$ is the origin
in $\mathbb R^3$, $|A| = |D| = 1$, and the angle between the
segments $OA$ and $OD$ is equal to $\frac{\pi}{6}$.

The manifold
$Y$
will be constructed so that
$O \in \partial Y$,
and
$]O,A] \subset \inter Y$,
$]O,D] \subset \inter Y$.
Under these conditions,
$\rho_Y(O,A) = \rho_Y(O,D) =1$.
However, the boundary of
$Y$
will create ``obstacles'' between
$A$
and
$D$
such that the length of any curve joining
$A$
and
$D$
in
$\inter Y$
will be greater than
$\frac{12}{5}$
(this means the violation of the triangle inequality for~$\rho_Y$).

Consider a countable collection of mutually disjoint segments
$\{I^k_j\}_{j \in \mathbb N,\, k = 1,\dots,k_j}$ lying in the
interior of the triangle $6 \Delta AOD$ (which is obtained from
the original triangle $\Delta AOD$ by dilation with
coefficient~$6$) with the following properties:

$(*)$
every segment
$I^k_j =[x^k_j,y^k_j]$
lies on a ray starting at the origin,
$y^k_j = 11 x^k_j$,
and
$|x^k_j| = 2^{-j}$;

$(**)$
For any curve
$\gamma$
with ends
$A$,
$D$
whose interior points lie in the interior of the triangle
$4 \Delta AOD$
and belong to no~segment
$I^k_j$,
the estimate
$l(\gamma) \ge 6$
holds.

The existence of such a family of segments is certain: they must be situated
chequerwise so that any curve disjoint from them be sawtooth, with the total length of
its ``teeth'' greater than $6$ (it can clearly be made greater than any prescribed positive
number). However, below we exactly describe the construction.

It is easy to include the above-indicated family of segments in the boundary
$\partial Y$
of~$Y$.
Thus, it creates a desired ``obstacle'' to joining
$A$
and
$D$
in the plane of
$\Delta AOD$.
But it makes no obstacle to joining
$A$
and
$D$
in the space. The simplest way to create such a~space obstacle
is as follows: Rotate each segment
$I^k_j$
along a spiral around the axis
$OA$.
Make the number of coils so large that the length of
this spiral be large and its pitch (i.e., the distance
between the origin and the end of a coil) be sufficiently small.
Then the set
$S^k_j$
obtained as the result of the rotation of the segment
$I^k_j$
is diffeomorphic to a plane rectangle, and it lies in a small neighborhood of
the cone of revolution with axis
$AO$
containing the segment
$I^k_j$.
The last circumstance guarantees that the sets
$S^k_j$
are disjoint as before, and so (as above) it is easy
to include them in the boundary
$\partial Y$
but, due to the properties of the
$I^k_j$'s
and a~large number of coils of the spirals
$S^k_j$,
any curve joining
$A,D$
and disjoint from each
$S^k_j$
has length
$\ge \frac{12}{5}$.

We turn to an exact description of the constructions used. First describe the construction
of the family of segments
$I^k_j$.
They are chosen on the basis of the following observation:

Let $\gamma : [0,1] \to 4 \Delta AOD$ be any curve with ends
$\gamma (0) = A$, $\gamma (1) = D$ whose interior points lie in
the interior of the triangle $4 \Delta AOD$. For $j \in \mathbb
N$, put $R_j = \{x \in 4 \Delta AOD : |x| \in [8 \cdot 2^{-j}, 4
\cdot 2^{-j}]\}$. It is clear that
$$
4 \Delta AOD \setminus \{O\} = \cup_{j \in \mathbb N} R_j.
$$
Introduce the polar system of coordinates on the plane of the
triangle $\Delta AOD$ with center $O$ such that the coordinates of
the points $A,D$ are $r = 1$, $\varphi = 0$ and $r = 1$, $\varphi
= \frac{\pi}{6}$, respectively. Given a point $x \in 6 \Delta
AOD$, let $\varphi_x$ be the angular coordinate of $x$ in
$[0,\frac{\pi}{6}]$. Let $\Phi_j = \{\varphi_{\gamma(t)} :
\gamma(t) \in R_j\}$. Obviously, there is $j_0 \in \mathbb N$ such
that
\begin{equation}\label{2}
\qquad\qquad\qquad\qquad\qquad\mathcal H^1(\Phi_{j_0}) \ge
2^{-j_0} \frac{\pi}{6},
\end{equation}
where $\mathcal H^1$ is the Hausdorff $1$-measure. This means
that, while in the layer $R_{j_0}$, the curve $\gamma$ covers the
angular distance $\ge 2^{-j_0} \frac{\pi}{6}$. The segments
$I^k_j$ must be chosen such that~(\ref{2}) together with the
condition
$$
\gamma(t) \cap I^k_j = \varnothing \quad \forall t \in [0,1]\,\,
\forall j \in \mathbb N\,\, \forall k \in \{1,\dots,k_j\}
$$
give the desired estimate
$l(\gamma) \ge 6$.
To this end, it suffices to take
$k_j = [(2 \pi)^j]$
(the integral part of
$(2 \pi)^j$)
and
$$
I^k_j = \{x \in 6 \Delta AOD : \varphi_x = k (2 \pi)^{-j} \frac{\pi}{6},\,
|x| \in [11 \cdot 2^{-j},2^{-j}]\},
$$
$k = 1,\dots,k_j$. Indeed, under this choice of the $I^k_j$'s,
estimate~(\ref{2}) implies that $\gamma$ must intersect at least
$(2 \pi)^{j_0} 2^{-j_0} = \pi^{j_0} > 3^{j_0}$ of the figures
$$
U_k = \{x \in R_{j_0} : \varphi_x \in (k(2 \pi)^{-j_0} \frac{\pi}{6},
(k + 1) (2 \pi)^{-j_0} \frac{\pi}{6})\}.
$$
Since these figures are separated by the segments $I^k_{j_0}$ in
the layer $R_{j_0}$, the curve $\gamma$ must be disjoint from them
each time in passing from one figure to another. The number of
these passages must be at least $3^{j_0} - 1$, and a~fragment of
$\gamma$ of length at least $2 \cdot 3 \cdot 2^{-j_0}$ is required
for each passage (because the ends of the segments $I^k_{j_0}$ go
beyond the boundary of the layer $R_{j_0}$ containing the figures
$U_k$ at distance $3 \cdot 2^{-j_0}$). Thus, for all these
passages, a section of $\gamma$ is spent of~length at least
$$
6 \cdot 2^{-j_0} (3^{j_0} -1) \ge 6.
$$
Hence, the construction of the segments $I^k_j$ with the
properties~($*$)--($**$) is finished.

Let us now describe the construction of the above-mentioned space
spirals.

For $x \in \mathbb R^3$, denote by $\Pi_x$ the plane that passes
through $x$ and is perpendicular to the segment $OA$. On
$\Pi_{x^k_j}$, consider the polar coordinates
$(\rho, \psi)$ with origin at the point of intersection
$\Pi_{x^k_j} \cap [O,A]$ (in this system, the point $x^k_j$ has
coordinates $\rho =\rho^k_j$, $\psi = 0$). Suppose that a~point
$x(\psi) \in \Pi_{x^k_j}$ moves along an~Archimedean spiral,
namely, the polar coordinates of~$x(\psi)$ are
$\rho(\psi) = \rho^k_j - \varepsilon_j \psi$, $\psi
\in [0,2 \pi M_j]$, where $\varepsilon_j$ is a small parameter to
be specified below, and $M_j \in \mathbb N$ is chosen so large
that the length of any curve passing between all coils of the
spiral is at least $10$.

Describe the choice of
$M_j$
more exactly. To this end, consider the points
$x(2 \pi)$,
$x(2 \pi (M_j - 1))$,
$x(2 \pi M_j)$,
which are the ends of the first, penultimate, and last coils of the spiral
respectively (with $x(0) = x^k_j$ taken as the starting point of the spiral).
Then
$M_j$
is chosen so large that the following condition hold:

$(*_1)$
{\it The length of any curve on the plane
$\Pi_{x^k_j},$
joining the segments
$[x^k_j,x(2 \pi)]$
and
$[x(2 \pi (M_j - 1)),x(2 \pi M_j)]$
and disjoint from the spiral
$\{x(\psi) : \psi \in [0,2 \pi M_j]\},$
is at least
$10$.}

Figuratively speaking, the constructed spiral bounds a~``labyrinth'', the mentioned
segments are the entrance to and the exit from this labyrinth, and thus any path through
the labyrinth has length
$\ge 10$.

Now, start rotating the entire segment $I^k_j$ in space along the
above-mentioned spiral, i.e., assume that $I^k_j(\psi) = \{y =
\lambda x(\psi) : \lambda \in[1,11]\}$. Thus, the segment
$I^k_j(\psi)$ lies on the ray joining $O$ with $x(\psi)$ and has
the same length as the original segment $I^k_j = I^k_j(0)$. Define
the surface $S^k_j = \cup_{\psi \in [0, 2 \pi M_j]} I^k_j(\psi)$.
This surface is diffeomorphic to a plane rectangle (strip). Taking
$\varepsilon_j > 0$ sufficiently small, we may assume without loss
of generality that $2 \pi M_j \varepsilon_j$ is substantially less
than $\rho^k_j$; moreover, that the surfaces $S^k_j$ are mutually
disjoint (obviously, the smallness of $\varepsilon_j$ does not
affect property~$(*_1)$ which in fact depends on $M_j$).

Denote by $y(\psi)=11x(\psi)$ the second end of the segment
$I^k_j(\psi)$. Consider the trapezium $P^k_j$ with vertices
$y^k_j$, $x^k_j$, $x(2 \pi M_j)$, $y(2 \pi M_j)$ and sides
$I^k_j$, $I^k_j(2 \pi M_j)$, $[x^k_j,x(2 \pi M_j)]$, and
$[y^k_j,y(2 \pi M_j)]$ (the last two sides are parallel since they
are perpendicular to the segment $AO$). By construction, $P^k_j$
lies on the plane $AOD$; moreover, taking $\varepsilon_j$
sufficiently small, we can obtain the situation where the
trapeziums $P^k_j$ are mutually disjoint (since $P^k_j \to I^k_j$
under fixed $M_j$ and $\varepsilon_j \to 0$). Take an arbitrary
triangle whose vertices lie on $P^k_j$ and such that one of these
vertices is also a vertex at an~acute angle in~$P^k_j$. By
construction, this acute angle is at least $\frac{\pi}{2} - \angle
AOD = \frac{\pi}{3}$. Therefore, the ratio of the side of
the~triangle lying inside the trapezium $P^k_j$ to the sum of the
other two sides (lying on the corresponding sides of $P^k_j$) is
at least $\frac{1}{2} \sin \frac{\pi}{3} > \frac{2}{5}$. If we
consider the same ratio for the case of a~triangle with a vertex
at an~obtuse angle of $P^k_j$ then it is greater than
$\frac{1}{2}$. Thus, we have the following property:

$(*_2)$ {\it For arbitrary triangle whose vertices lie on the
trapezium $P^k_j$ and one of these vertices is also a~vertex in
$P^k_j$, the sum of lengths of the sides situated on the
corresponding sides of $P^k_j$ is less than $\frac{5}{2}$ of the
length of the third side} ({\it lying inside} $P^k_j$).

Let a point $x$ lie inside the cone $K$ formed by the rotation of
the angle $\angle AOD$ around the ray $OA$. Denote by $\Proj x$
the point of the angle $\angle AOD$ which is the image of $x$
under this rotation. Finally, let $K_{4 \Delta AOD}$ stand for the
corresponding truncated cone obtained by the rotation of the
triangle $4 \Delta AOD$, i.e., $K_{4 \Delta AOD} = \{x \in K :
\Proj x \in 4 \Delta AOD\}$.

The key ingredient in the proof of our theorem is the following
assertion:

$(*_3)$ {\it For arbitrary space curve $\gamma$ of length less
than $10$ joining the points $A$ and $D$, contained in the
truncated cone $K_{4 \Delta AOD} \setminus \{O\}$, and disjoint
from each strip $S^k_j$, there exists a plane curve
$\tilde{\gamma}$ contained in the triangle $4 \Delta AOD \setminus
\{O\},$ that joins $A$ and $D$, is disjoint from all segments
$I^k_j$ and such that the length of $\tilde{\gamma}$ is  less than
 $\frac{5}{2}$ of the length of} $\Proj \gamma$.

Prove~$(*_3)$. Suppose that its hypotheses are fulfilled. In
particular, assume that the inclusion $\Proj \gamma \subset 4
\Delta AOD \setminus \{O\}$ holds. We need to modify $\Proj
\gamma$ so that the new curve be contained in the same set but be
disjoint from each of the $I^k_j$'s. The construction splits into
several steps.

{\bf Step~1.} If $\Proj \gamma$ intersects a segment $I^k_j$ then
it necessarily intersects also at least one of the shorter sides
of $P^k_j$.

Recall that, by construction, $P^k_j = \Proj S^k_j$; moreover,
$\gamma$ intersects no spiral strip $S^k_j$. If $\Proj \gamma$
intersected $P^k_j$ without intersecting its shorter sides then
$\gamma$ would pass through all coils of the corresponding spiral.
Then, by~$(*_1)$, the length of the corresponding fragment of
$\gamma$ would be $\ge 10$ in contradiction to our assumptions.
Thus, the assertion of step~1 is proved.

{\bf Step~2.} Denote by $\gamma_{P^k_j}$ the fragment of the plane
curve $\Proj \gamma$ beginning at the first point of its entrance
into the trapezium $P^k_j$ to the point of its exit from $P^k_j$
(i.e., to its last intersection point with $P^k_j$). Then this
fragment $\gamma_{P^k_j}$ can be deformed without changing the
first and the last points so that the corresponding fragment of
the new curve lie entirely on the union of the sides of $P^k_j$;
moreover, its length is at most $\frac{5}{2}$ of the length of
$\gamma_{P^k_j}$.

The assertion of step~2 immediately follows from the assertions of step~1 and
$(*_2)$.

The assertion of step~2 in turn directly implies the desired assertion
$(*_3)$.
The proof of
$(*_3)$
is finished.

Now, we are ready to pass to the final part of the proof of Theorem~2.

$(*_4)$ {\it The length of any space curve $\gamma \subset \mathbb
R^3 \setminus \{O\}$ joining $A$ and $D$ and disjoint from each
strip $S^k_j$ is at least} $\frac{12}{5}$.

Prove the last assertion. Without loss of generality, we may also
assume that all interior points of $\gamma$ are inside the cone
$K$ (otherwise the initial curve can be modified without any
increase of its length so that it have property~$(*_4)$). If
$\gamma$ is not included in the truncated cone $K_{4 \Delta AOD}
\setminus \{O\}$ then $\Proj \gamma$ intersects the segment
$[4A,4D]$; consequently, the length of $\gamma$ is at least $2(4
\sin \angle OAD - 1) = 2(4 \sin \frac{\pi}{3} - 1) = 2 (2 \sqrt 3
- 1) > 4$, and the desired estimate is fulfilled. Similarly, if
the length of $\gamma$ is at least $10$ then the desired estimate
is fulfilled automatically, and there is nothing to prove. Hence,
we may further assume without loss of generality that $\gamma$ is
included in the truncated cone $K_{4 \Delta AOD} \setminus \{O\}$
and its length is less than $10$. Then, by~$(*_3)$, there is a
plane curve $\tilde \gamma$ contained in the triangle $4 \Delta
AOD \setminus \{O\}$, joining the points $A$ and $D$, disjoint
from each segment $I^k_j$, and such that the length of $\tilde
\gamma$ is at most $\frac{5}{2}$ of the length of $\Proj \gamma$.
By property~$(**)$ of the family of segments $I^k_j$, the length
of $\tilde \gamma$ is at least $6$. Consequently, the length of
$\Proj \gamma$ is at least $\frac{12}{5}$, which implies the
desired estimate. Assertion~$(*_4)$ is proved.

The just-proven property $(*_4)$ of the constructed objects
implies Theorem~2. Indeed, since the strips $S^k_j$ are mutually
disjoint and, outside every neighborhood of the origin $O$, there
are only finitely many of these strips, it is easy to construct a
$C^0$-manifold $Y \subset \mathbb R^3$ that is homeomorphic to a
closed ball (i.e., $\partial Y$ is homeomorphic to a
two-dimensional sphere) and has the following properties:

(I)
$O \in \partial Y$,
$[A,O[ \cup [D,O[ \subset \inter Y$;

(II) for every point $y \in (\partial Y) \setminus \{O\}$, there
exists a neighborhood $U(y)$ such that $U(y) \cap \partial Y$ is
$C^1$-diffeomorphic to the plane square $[0,1]^2$;

(III) $S^k_j \subset \partial Y$ for all $j \in \mathbb N,\, k =
1,\dots,k_j$.

The construction of $Y$ with properties~(I)--(III) can be carried
out, for example, as follows: As the surface of the zeroth step,
take a~sphere containing $O$ and such that $A$ and $D$ are inside
the sphere. On the $j$th step, a small neighborhood of the point
$O$ of our surface is smoothly deformed so that the modified
surface is still smooth, homeomorphic to a~sphere, and contains
all strips $S^k_j$, $k = 1,\dots,k_j$. Besides, we make sure that,
at the each step, the so-obtained surface be disjoint from the
half-intervals $[A,O[$ and $[D,O[$, and, as above, contain all
strips $S^k_i$, $i \le j$, already included therein. Since the
neighborhood we are deforming contracts to the point $O$ as $j \to
\infty$, the so-constructed sequence of surfaces converges (for
example, in the Hausdorff metric) to a limit surface which is the
boundary of a $C^0$-manifold $Y$ with properties~(I)--(III).

Property~(I) guarantees that
$\rho_Y(A,O) = \rho_Y(A,D) = 1$
and
$\rho_Y(O,x) \le 1 + \rho_Y(A,x)$
for all
$x \in Y$.
Property~(II) implies the estimate
$\rho_Y(x,y) < \infty$
for all
$x,y \in Y \setminus \{O\}$,
which, granted the previous estimate, yields
$\rho_Y(x,y) < \infty$
for all
$x,y \in Y$.
However, property~(III) and the assertion~$(*_4)$
imply that
$\rho_Y(A,D) \ge \frac{12}{5} > 2 = \rho_Y (A,O) + \rho_Y(A,D)$.
Theorem~2 is proved.

In the case where
$\rho_Y$
is a metric (the dimension
$n$
($\ge 2$)
is arbitrary), the question of the existence of geodesics is solved
in the following assertion, which implies that
$\rho_Y$
is an {\it intrinsic metric} (see, for example, \S 6
from~\cite{al1}). 

{\bf Theorem~3.}
{\it Assume that
$\rho_Y$
is a finite function and is a metric on~$Y$.
Then any two points
$x,y \in Y$
can be joined
in~$Y$
by a~shortest curve
$\gamma : [0,L] \to Y$
in the metric
$\rho_Y;$
i.e.{\rm,}
$\gamma(0) = x,$
$\gamma(L) = y$,
and
\begin{equation}
\qquad\qquad\rho_Y(\gamma(s),\gamma(t)) = t - s \quad \forall s,t
\in [0,L], \quad s < t. \label{3} \end{equation} }

{\sc Proof.} Fix a pair of distinct points $x,y \in Y$ and put $L
= \rho_Y (x,y)$. Now, take a~sequence of paths $\gamma_j : [0,L]
\to Y$ such that $\gamma_j(0) = x_j$, $\gamma_j(L) = y_j$, $x_j
\to x$, $y_j \to y$, and $l(\gamma_j) \to L$ as $j \to \infty$.
Without loss of generality, we may also assume that the
parametrizations of the curves $\gamma_j$ are their natural
parametrizations up to a~factor (tending to~$1$) and the mappings
$\gamma_j$ converge uniformly to a~mapping $\gamma : [0,L] \to Y$
with $\gamma (0) = x$, $\gamma (L) = y$. By these assumptions,
\begin{equation}
\qquad\qquad\lim_{j \to \infty} l(\gamma_j|_{[s,t]}) = t - s \quad
\forall s,t \in [0,L], \quad s < t. \label{4} \end{equation}

Take an arbitrary pair of numbers $s,t \in [0,L]$, $s < t$. By
construction, we have the convergence $\gamma_j(s) \in \inter Y
\to \gamma(s)$, $\gamma_j(t) \in \inter Y \to \gamma(t)$ as $j \to
\infty$. From here and the definition of the metric
$\rho_Y(\cdot,\cdot)$ it follows that
$$
\rho_Y(\gamma(s),\gamma(t)) \le \lim_{j \to \infty}
l(\gamma_j|_{[s,t]}).
$$
By~(\ref{4}), \begin{equation}
\quad\qquad\qquad\rho_Y(\gamma(s),\gamma(t)) \le t - s \quad
\forall s,t \in [0,L],\,\, s < t. \label{5}
\end{equation} Prove that~(\ref{5}) is indeed an equality. Assume that
$$
\rho_Y(\gamma(s'),\gamma(t')) < t' -s'
$$
for some $s',t' \in [0,L]$, $s' < t'$. Then, applying the triangle
inequality and then~(\ref{5}), we infer
$$
\rho_Y(x,y) \le \rho_Y(x,\gamma(s')) +
\rho_Y(\gamma(s'),\gamma(t')) + \rho_Y(\gamma(t'),y)
< s' + (t' - s') + (L - t') = L,
$$
which contradicts the initial equality $\rho_Y(x,y) = L$. The
so-obtained contradiction completes the proof of
identity~(\ref{3}).

{\bf Remark.} Identity~(\ref{3}) means that the curve $\gamma$ of
Theorem~3 is a~geodesic in the metric $\rho_Y$, i.e., the length
of its fragment between points $\gamma(s)$, $\gamma(t)$ calculated
in~$\rho_Y$ is equal to $\rho_Y(\gamma(s),\gamma(t)) = t - s$.
Nevertheless, if we compute the length of the above-mentioned
fragment of the curve in the initial Riemannian metric then this
length need not coincide with $t-s$; only the easily verifiable
estimate $l(\gamma|_{[s,t]}) \le t - s$ holds (see~(\ref{4})\,).
In the general case, the equality $l(\gamma|_{[s,t]}) = t - s$ can
only be guaranteed if $n = 2$ (if $n \ge 3$ then the corresponding
counterexample is constructed by analogy with the counterexample
in the proof of Theorem~2, see above). In particular, though, by
Theorem~3, the metric $\rho_Y$ is always intrinsic in the sense of
the definitions in~\cite[\S 6]{al1}, 
the space $(Y,\rho_Y)$ may fail to
be {\it a space with intrinsic metric} in the sense of [ibid].

\vskip3mm

{\small
\section*{Acknowledgements}

The authors were partially supported by the Russian Foundation for
Basic Research (Grant~11-01-00819-a), the Interdisciplinary
Project of the Siberian and Far-Eastern Divisions of the Russian
Academy of Sciences (2012-2014 no.~56), the State Maintenance
Program for the Leading Scientific Schools of the Russian
Federation (Grant~NSh-921.2012.1), and the Exchange Program
between the Russian and Polish Academies of Sciences (Project
2011--2013). }

{\small


}

\end{document}